\documentclass{ifacconf}

\usepackage{graphicx}      
\usepackage{natbib}        
\usepackage{amsmath,amssymb,bm}
\usepackage{array}
\usepackage{pst-all} 
\usepackage{relsize} 
\def\qed{
\vbox{\hrule\hbox{\vrule\hbox to 5pt{\vbox to 8pt{\vfil}\hfil}\vrule}\hrule}}
\def\qedsmall{
\vbox{\hrule\hbox{\vrule\hbox to 4pt{\vbox to 4pt{\vfil}\hfil}\vrule}\hrule}}
\def\endproof{\unskip \nobreak \hskip0pt plus 1fill \qquad \qed}
\def\endrem{\unskip \nobreak \hskip0pt plus 1fill \qquad \qedsmall}
\def\es{\enspace}

\def\9{\"u9}

\def\R{\mathbb{R}}
\def\S{\mathbb{S}}

\def\N{\mathbb{N}}

\def\KK{{\cal K}}

\def\LL{{\cal L}}

\def\eps{\varepsilon}

\def\beq{\begin{equation}}
\def\eeq{\end{equation}}
\def\bea{\begin{eqnarray}}
\def\eea{\end{eqnarray}}
\def\beaa{\begin{eqnarray*}}
\def\eeaa{\end{eqnarray*}}

\newtheorem{tsatz}{Theorem}[section]\newenvironment{definition}{\begin{tdefi} \em}{\em \endrem \end{tdefi}}

\newtheorem{tdefi}[tsatz]{Definition}
\newtheorem{tbem}[tsatz]{Remark}

\newtheorem{tann}[tsatz]{Assumption}

\newenvironment{theorem}{\begin{tsatz} \em}{\em \end{tsatz}}
\newenvironment{remark}{\begin{tbem} \em}{\em \endrem \end{tbem}}
\newenvironment{assumption}{\begin{tann} \em}{\em \endrem \end{tann}}

\input epsf
\def\beq{\begin{equation}}
\def\eeq{\end{equation}}
\def\bea{\begin{eqnarray}}
\def\eea{\end{eqnarray}}
\def\beaa{\begin{eqnarray*}}
\def\eeaa{\end{eqnarray*}}

\begin{document}
\begin{frontmatter}

\title{Robustness of Prediction Based Delay Compensation for Nonlinear Systems}

\author[First]{Rolf Findeisen}
\author[Second]{Lars Gr\"une}
\author[Third]{J\"urgen Pannek}
\author[First]{Paolo Varutti}

\address[First]{  Otto--von--Guericke--Universit{\"a}t Magdeburg,
  Institute for Automation Engineering, Chair for Systems Theory and
  Automatic Control, \{rolf.findeisen,paolo.varutti\}@ovgu.de}  
\address[Second]{University of Bayreuth, Institute of Mathematics,
  Chair of Applied Math., lars.gruene@uni-bayreuth.de} 
\address[Third]{Curtin University of Technology, Perth, Australia,
  juergen.pannek@googlemail.com}
\begin{abstract}                
Control of systems where the information between the controller,
actuator, and sensor can be lost or delayed can be challenging with
respect to stability and performance. One way to 
overcome the resulting problems  is the use of prediction based
compensation  schemes. Instead of a single input, a sequence of (predicted) future 
controls is submitted and implemented at the actuator. If suitable,
so-called prediction consistent compensation and control schemes, such
as certain predictive control approaches, are used,
stability of the closed loop in the presence of delays and packet losses can be
guaranteed. In this paper, we show that control schemes employing prediction
based delay compensation approaches do posses inherent
robustness properties. Specifically, if the nominal closed loop 
system without delay compensation is ISS with respect to perturbation and
measurement errors, then the closed loop system employing prediction based delay
compensation techniques is robustly stable. We analyze the
influence of the prediction horizon on the robustness gains and illustrate the
results in simulation.
\end{abstract}


\end{frontmatter}
\section{Introduction}\label{sec:introduction}
In many of today's control systems, delays and information losses
between the controller, actuator, and sensor are often
unavoidable. Such delays and dropouts must be accounted for during the
controller 
design and the closed loop system analysis to avoid instability or performance
decay. There are many causes for
delays and information losses:  often controller,
sensors, and actuators are connected via a communication
network. These kinds of systems are typically denoted as networked control
systems, see \cite{Hespanha2007}. For networked control systems delays and information losses 
might be  due to network overload, long communication
distances, routing, or hardware failure. 
Other sources might be long computation times, e.g. due to
image processing, or solution of complex optimization problem,  
as in the case of predictive control. 
In other applications, delays and potential
losses might be inherent to the problem under investigation. 
For instance, in some processes it is necessary to ``charge'' batteries
before being able to collect a measurement or human interaction could be required to take measurements or
apply new inputs. \\
By now, a series of approaches for the control of systems that are
subject to delays and information losses exist. 
In particular, ideas based on Model Predictive Control (MPC)
have demonstrated to be effective in dealing with both delays and
information losses, see
e.g. \cite{Findeisen2004,Findeisen2006a,Findeisen2009,Varutti2009d,Grune2009c,Polushin2008,Lunze2010,Bemp98}. 
Many of these approaches are based on the idea of compensating the unknown
delays by sending not only one control action to the actuator,
rather a complete input sequence (discrete time systems) or an input
signal (continuous time systems), are submitted, where the signals
sent are time-stamped.  
The actuator itself can then continue applying the old input
until new data arrive or the input has been implemented.
While suitable design of such a scheme often leads to nominal closed loop
stability, only minor results with respect to robust stability are
available.\\
In this work, we  establish robust stability properties of  so called
prediction consistent delay compensation schemes, see \cite{Findeisen2009} and \cite{Varutti2009b}
for the continuous time formulations and \cite{Grune2009c} for the
discrete time case. Prediction consistent delay compensation  schemes counteract the delays by submitting
a complete consistent input trajectory (in the sense of a predicted behavior of the closed loop) to the actuator.
Approaches following similar ideas, which, however, are only able to handle delays either
on the sensor or actuator side, have been for example introduced in
\cite{Bemp98,Polushin2008} (see also \cite{Grune2009c} for a comparison).\\
Specifically, we establish that prediction consistent delay
compensation schemes for discrete time nonlinear systems do admit,
under certain conditions, inherent robustness. Precisely, if
the nominal closed loop system without delays
is input-to-state stable (ISS), cf. \cite{Sontag2000},
then the closed loop system subject to delays and utilizing a prediction consistent delay
compensation approach is robustly stable. 
Additionally, we explicitly analyze the influence of the prediction horizon on
the robustness gains and illustrate the results by a numerical
simulation. The derived results significantly expand the applicability
of prediction consistent delay compensation approaches, 
since they establish that these methods
are also well suited for the robust case. It is important to stress
that ISS results for predictive control methods are well
established by now, refer to \cite{magni07a,limon09a}.
Results with respect to ISS, delays, and predictive control 
approaches are, however, very limited. Exceptions
are the results presented in \cite{zavala09}, which, however, do not
apply a compensation approach and thus can only derive ISS properties
with respect to small delays. Similar results with respect to
practical stability subject to delays have been presented in
\cite{Findeisen2006a}. Furthermore, recently stochastic stability 
properties of  
predictive control approaches over unreliable networks have been
derived in \cite{quevedo2009}. These results, however, only hold for delays
and losses of information on the actuator side.  
\section{Prediction Consistency: Setting and Nominal Results}\label{sec:setting-notation}
In this paper, we consider  discrete time nonlinear systems
\beq x(n+1) = f(x(n),u(n),w(n)), \label{sys}\eeq
where $x \!\in\! X\! \subseteq \!\R^n$, $u \!\in\! U \!\subseteq \!\R^m$, and $w \!\in\! W$ represent respectively the
state, the input and the disturbance/perturbation acting on the system, taken from
the sets $X$, $U$, and $W$. In the following, for any  time $n>n_0$,  
we denote with $x(n,n_0,x_0,u,w)$ the solution of  \eqref{sys} with initial time $n_0$, initial value $x_{n_0}$,
input sequences $u$ and perturbation $w$, obtained by iterating
\eqref{sys}  for $n_0,\ldots, n$.
\begin{remark}
  The disturbance $w$ can account for various uncertainties
  such as measurement uncertainties, disturbances on the
  input side, or model uncertainties.
\end{remark}

\begin{figure}[htb]
\centering
\includegraphics[width=0.45\textwidth]{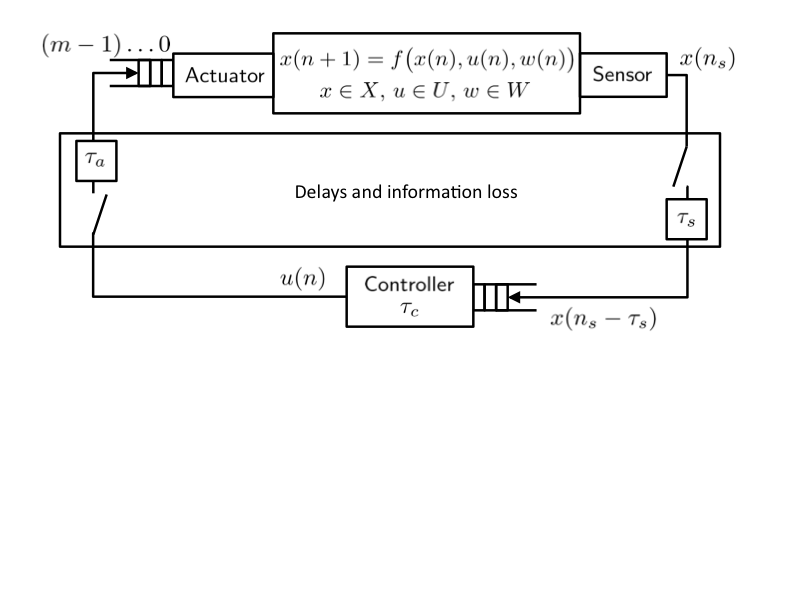}\vspace*{-2.7cm}
\caption{Sketch of the problem under investigation.\vspace{-.6cm}}\label{fig:system}
\end{figure}

{\bf Prediction consistent compensation:}\\
We consider that the controller interacts with the sensors and
actuators  as shown in Figure~\ref{fig:system}. The implementation of the controller  follows
the ideas presented in \cite{Grune2009c} for discrete time systems and
in \cite{Varutti2009b} for the continuous time case. 
A key component of all these approaches is a prediction algorithm,
which makes use of a model of the system to deal with delays and
information losses by forward predicting the (closed loop) plant and
generating an input sequence which is communicated to the actuator. 

For the prediction, a model of the plant is required. Our scheme relies
on an approximating map $\tilde f(\cdot)$ of the  nominal unperturbed
 system map $f(\cdot)$ in the sense that   
\beq
\tilde f(x,u) \approx f(x,u,0)
\footnote{If desired and available, an estimate $\hat w \approx w$
 could be used in the prediction. However, this will make the
 subsequent analysis somewhat more involved, in particular for the
 optimization based MPC case. Furthermore, in practice suitable
 estimates $\widehat w \approx w$ will only be available in
 exceptional cases, which might lead to computational demanding
 $min-max$ or set-based approaches.}.
\eeq
The detailed approximation properties of $\tilde f(\cdot)$ are
elaborated in Assumption \ref{errass}(i), below. For instance, if $f(\cdot)$ is a
discrete time model of a zero order hold sampled continuous time
control system, then $\tilde f(\cdot)$ might be chosen as the numerical
solution of the underlying differential equation over one sampling
period with constant control value. Analogously to the plant, 
$\tilde x(n,n_0,\tilde x_{n_0},\tilde u)$ denotes the predicted solutions
obtained by iterating  
\beq \tilde x(n+1) = \tilde f(\tilde x(n),\tilde u(n))
\label{asys}\eeq
from $n_0$ to $n$ starting with initial
value $\tilde x_{n_0}$ using $\tilde u$. \\
This prediction map is used by the controller to ``forecast''
the future system states based on past measurements. To infer
properties of the closed loop from these forecast predictions, the
real system state should coincide, or at least be close to the
predicted state. A key requirement for the prediction $\tilde x$ to
assume the same value as the real system state $x$ is that the control
sequence $\tilde u$ used for the prediction coincides 
with the control sequence $u$ applied at the plant. In other words, we require the
compensation algorithm to be \emph{prediction consistent}, or, more
formally:
\begin{definition}\label{sec:pred-cons-sett}(Prediction Consistency):\\
 A prediction based delay compensation control scheme is
called {\em prediction consistent} if
at each time $n\in\N_0$ the identity 
$\tilde u(k) = u(k)$ holds for all $k\in \{0,\ldots,n\}$.
\end{definition}
Note that the schemes proposed in
\cite{Bemp98,Grune2009c,Polushin2008} are prediction consistent in the
sense of this definition.  \\
Before describing in detail the scheme used, we
introduce the following quantities and assumptions on the delays and
information losses. 

We refer with $n_s$, $n_c$, and $n_a$ to the time at the sensor,
controller and actuator side respectively;
$\tau_s$, and $\tau_a$ indicate the delays in the measurement
and actuation channel, while $\tau_c$ represents potential
computational delays; all delays are assumed to be
bounded\footnote{This assumption can be easily relaxed by defining
  countdowns after which the exchanged information is 
  considered as lost.}. Note, however, that the delays are not assumed
to be constant, they can vary with time. Both measurements and input sequences can be
lost 
and all the components --- controller, sensor,
and actuator --- are supposed to have synchronized clocks. All
measurements are time-stamped with the time instant $n_s$ from which
they are collected, whereas the input is time-stamped
with the value $n$ in which it is supposed  
to be applied to the system.
The disordered arrival of packets is solved by taking the one with the
most recent time-stamp. Additionally, the actuator is equipped with a
buffer of  length $m>1$.\\

To counteract and compensate delays and information losses, we propose
the following scheme:
\begin{enumerate} 
\item[(i)] At each  time $n_c$ the controller computes the
predicted state $\tilde x(n)=\tilde x(n,n_s, x(n_s),\tilde u)$ at time
$n$ (specified in Step (ii)) from the delayed measurement
$x(n_s)$ taken at time $n_s$ using the prediction
control input $\tilde u$. 
\item[(ii)] Based on the prediction $\tilde x(n)$, and the buffer
  length $m$, the controller computes  an input  sequence 
\beq \mu(\tilde x(n),0),\ldots,\mu(\tilde x(n),m-1) \label{museq}\eeq
which is sent to the actuator. 
This allows to have ``backup'' control values whenever transmission
failures in the input or long delays occur.The prediction time $n$ is chosen such that the input 
sequence reaches the actuator in time. More precisely, for the
prediction performed at time $n_c$, using an upper
bound $\tau^{\max}\ge \tau_a+\tau_c$ we predict the future state for
time $n=n_c + \tau^{\max} > n_c$ from the most recent available measurement
taken at time $n_s \le n_c$. The resulting input sequence is
time-stamped with  $n$ and sent to the actuator. The value
$\tau(n):=n-n_s$ denotes the length of the resulting prediction
interval. 
\item[(iii)] The actuator buffers the input sequences and uses the most recent
 available one in its buffer to determine the control value applied to
 the plant.  
\end{enumerate}
Different approaches to generate prediction consistent input sequences
can be employed,  as long as they satisfy the prediction consistency
condition stated in Definition \ref{sec:pred-cons-sett}.
Without loss of generality, to generate the control sequence \eqref{museq}
of Step (ii) we consider two different approaches: generation by forward
predicting a stabilizing static feedback, and by using predictive
control schemes. Notice that it is possible to compute these
sequences in many other ways, 
cf.\ \cite{Grune2009c,Varutti2009b,Polushin2008}. All these
approaches can be analyzed with the framework proposed in this
paper.\\[1mm]
{\bf Prediction consistent input generators:}
\begin{itemize} 
\item[(iia)] {\it Static state feedback:} the input sequence can be
  based on a  static\footnote{The controller can in
    principle also be dynamic, which is avoided here for simplicity of
  presentation.} state feedback controller $K:X\to
  U$. In this case, from the predicted state $\tilde x(n)\in X$ we
  inductively compute
 \beq 
 \begin{array}{c}
 \tilde x(p+1) = \tilde f(\tilde x(p),K(\tilde x(p))),\\
 p=n,\ldots,n+m-2, \label{intpredict_a}
 \end{array}
 \eeq
 and set 
 \beq
 \mu(\tilde x(n),q) := K(\tilde x(q)), \quad
 q=0,\ldots,m-1. \label{Keval}
 \eeq 
 For notation simplicity, although this is in general not required, we assume that
 the predictor \eqref{intpredict_a} coincides with the one specified in Step (i).
 If we extend the prediction control sequence $\tilde u$ used in Step
 (i) for computing $\tilde x(n)$ by 
 setting $\tilde u(n+q)=K(\tilde x(q))$, $q=0,\ldots,m-1$, then from
 \eqref{intpredict_a},
\beq 
\begin{array}{c}
\tilde x(p) = \tilde x(p,n_s,x(n_s),\tilde u),\\
 p=n,\ldots,n+m-2, \label{xpalt}
 \end{array}
 \eeq 
 where $n_s=n-\tau(n)$ is the measurement time from Step (i).
\item[(iib)] {\it predictive control:} If the controller is computed by
  a predictive control approach (MPC),  then for a given initial value
  $x_0$ the MPC algorithm generates a finite horizon  control
  sequence  $u_{x_0}(0), \ldots, u_{x_0}(m-1)$. These control values are
 determined via an internal prediction inside the optimization algorithm
 starting from $x_0=\tilde x(n)$. This results in the predicted  optimal trajectory 
 \beq
 \begin{array}{c}
  \tilde x(p) = \tilde x(p, n, \tilde x(n), u_{\tilde
   x(n)}),\\
   p=n+1,\ldots, n+m-2. \label{mpcpred}
   \end{array}
 \eeq 
 For simplicity of exposition 
 we assume that
 \eqref{asys} is used to compute the internal prediction.
While in a usual MPC scheme one would only use the first element of
 the optimal control sequence for feedback, in order to
 obtain the sequence \eqref{museq}, we set  
 \beq \mu(\tilde x(n),q) = u_{\tilde x(n)}(q), \quad
 q=0,\ldots,m-1. \label{mufromKp}\eeq 
 Again we can write $\tilde x(p)$ in the form \eqref{xpalt} if we
 appropriately extend the prediction control sequence used in Step
 (i). Here we need to use
 $\tilde u(n+q)=u_{\tilde x(n)}(q)$, $q=0,\ldots,m-1$.
\end{itemize}
Both control approaches have certain advantages and
disadvantages. Usage of a known static feedback law allows, in
general, fast generation of input sequences. However, it might be
difficult to take constraints or cost functions to be optimized into
account. Obtaining a prediction consistent input sequence by MPC might
be computationally challenging, however, it allows to directly
consider constraints. 

{\it Overall closed loop system:}\\
By $\sigma_i$, $i\in\N_0$, we denote in the following the times ---
numbered in increasing order --- at which the actuator switches to a
new control sequence, i.e., the times at which 
a control value $\mu(\tilde x(\sigma_i),0)$ is applied at the
actuator. Henceforth, we will refer to the times $\sigma_i$ as the {\em
  switching  times}. Using this, we can write the closed loop as    
\beq
  x(n+1) = f(x(n),
  \mu ( \tilde x(\sigma_i),n-\sigma_i ), w(n)), 
  \label{cl2}
\eeq
for all $n \in\{\sigma_i,\sigma_{i+1}-1\}$, where $\tilde x(\sigma_i)= \tilde x(\sigma_i, \sigma_i-\tau(\sigma_i), x(\sigma_i-\tau(\sigma_i)),\tilde u)$, and $\tau(n)= n-n_s$.
The prediction consistency condition can now be ensured by suitable
algorithms which enable controller and actuator to identify and
correct prediction inconsistencies by means of sending time-stamped
information. For more details see \cite{Grune2009c} and
\cite{Varutti2009b}. If prediction consistency holds,
then from \eqref{cl2} it follows that
\beq
u(n) = \mu(\tilde x(\sigma_i,
\sigma_i-\tau(\sigma_i),x(\sigma_i-\tau(\sigma_i)),\tilde 
u),n-\sigma_i) \label{ucl}
\eeq
for all $n \in\{\sigma_i,\sigma_{i+1}-1\}$ and all
$i=0,1,2,\ldots$. 

In order to simplify the analysis, we shift our ``time''  and number
the $\sigma_i$ such that $\sigma_0=0$.
The resulting closed loop trajectory for
$n\ge 0$ is uniquely determined by the value $x_0=x(0)$, the switching
times $\sigma_i$ and the delays $\tau(\sigma_i)$. We denote this
closed loop trajectory by $x_{cl}(n,x_0,\sigma_\cdot,\tau,w)$ and use
the brief notation $x_{cl}(n)$.

\begin{remark} (Open loop prediction versus closed loop)\\
The closed loop system \eqref{cl2} appears to depend only on the
predictions $\tilde x(\sigma_i,
\sigma_i-\tau(\sigma_i),x(\sigma_i-\tau(\sigma_i)),\tilde 
u))$ for the switching times $\sigma_i$. However, in the Steps (iia)
and (iib), above, also the 
predictions $\tilde x(p)$ for all remaining times $p\ne \sigma_i$ are
needed. Note that each $p$ appears in (iia) or (iib) for several
different $n$ and thus at different times/runs of these steps different
values $\tilde x(p)$ are computed for one and the same $p$. In (iia)
only one of them, more precisely the one corresponding to the maximal 
$n$ satisfying $n=\sigma_i$ for some $i$, is actually used in
order to compute the control value $u(n)$ applied to the
closed loop system. Hence, for each closed loop trajectory
$x_{cl}(\cdot)$ and each time $n\ge 0$, there is a unique prediction
$\tilde x(n)$ which is used either explicitly for $n=\sigma_i$ in
\eqref{cl2} or implicitly for $n=p$ in Step (iia) in order to
compute the value $u(n)$  eventually applied to the system.
We denote this predicted state by $\tilde x_{cl}(n)$.   \\
In (iib), since the optimizations are carried out over the whole
prediction horizon, for each computation all predicted values (also
those for $k>n$) affect the computed control sequence. Since
these future predictions are uniquely determined by
$\tilde x_{cl}(n)$ and
\eqref{asys}, we do not denote them explicitly.
\label{rem:allpredict}\end{remark}

In the next section we establish the main result, namely
robust stability and the influence of the prediction horizon on the
robustness gains of the proposed scheme.

\section{Stability and robustness}\label{sec:stab-robustn-gener}
The main idea of our analysis is to replace the  closed loop system 
\eqref{cl2} by a non-delayed system in which the prediction errors due to the delay effects
are captured as measurement errors. This fundamentally  differs
from other stability analysis methods for delayed systems in
which the delay is explicitly taken into account. While our method
may lead to  more conservative results, its main advantage is the fact that it is
applicable to general nonlinear systems under rather mild conditions.

{\it Bounding the influence of  prediction errors:}\\
To capture the prediction errors via measurement errors, we need the
following additional assumptions with respect to the estimates on the prediction
accuracy and the sensitivity of the solution with respect to  $w$.
\begin{assumption}(Prediction error and perturbations)\\
(i) The prediction error for the nominal system satisfies
\[ \| \tilde x(n,n_0,\tilde x_0,u) -
x(n,n_0,x_0,u,0)\| \le \eps(n-n_0,\|\tilde x_0-x_0\|)\]
for all $x_0,x_1\in X$ and all $u\in U$, where $\eps:(\R_0^+)^2\to\R_0^+$
is a monotonically increasing, continuous function in both arguments. \\
(ii) The influence of the perturbation can be bounded by
\[ \|x(n,n_0,x_0,u,0) - x(n,n_0,x_0,u,w)\| \le \eta(n-n_0,
\|w\|_\infty)\]
for all $x_0\in X$ and all $u\in U$; $\eta:(\R_0^+)^2\to\R_0^+$
is a continuous function which is monotonically increasing in its first
argument and satisfies   
$\eta(n,\cdot)\in\KK$ for all $n\ge 0$. In particular, $\eta(n,w)=0$
for all $n\in\N_0$ if $w\equiv 0$.
\label{errass}\end{assumption}
If required, explicit expressions for $\eps$ and $\eta$ can be derived from
suitable properties for $f(\cdot)$ and $\tilde f(\cdot)$. For example, if $\tilde f(\cdot)$ is a
numerical approximation (e.g. for a continuous time system) with convergence order  
$p\in\N$ which is Lipschitz in $x$ with Lipschitz constant
$L$, then a standard error estimation from numerical analysis
yields 
\beq \eps(k,r) = (e^{L k}-1)K h^p+e^{Lk}r, \label{epsex}\eeq
where $h$ is the (fixed) time step used in the numerical scheme and $K$ is a
suitable constant. Alternatively, a numerical step size controlled
scheme could be used. In this case the term $Kh^p$ in \eqref{epsex} is
replaced by a user specified desired accuracy $\hat \eps$. 

If $f$ is Lipschitz in $x$ with constant $L$ and satisfies 
\[ \|f(x,u,w)-f(x,u,0)\| \le \rho(\|w\|) \]
for some $\KK_\infty$-function $\rho$, then 
\beq \eta(k,r) = (e^{L k}-1)\rho(r)/L \label{etaex}\eeq
holds.
\begin{remark} (Open loop stable and unstable systems)\\
  Note that the exponential growth of the error terms in $k$ in
  \eqref{epsex} and \eqref{etaex} is a worst
case estimate which applies if the plant is open loop unstable. 
If we assume that the system to be controlled is open loop stable for
$w\equiv 0$ (e.g. if the plant is pre-stabilized by a feedback
controller situated at the plant), then one only has linear, not  exponential
 growth of the error terms in $k$.  
\label{rem:open-loop-stable}\end{remark}

{\it Non-delayed closed loop system:}\\
The auxiliary system we use for the  analysis is given by 
\beq x(n+1) =
f(x(n),\mu(x(\sigma_i)+v(\sigma_i)),n-\sigma_i),w(n))\label{cla}\eeq 
for $n \in\{\sigma_i,\sigma_{i+1}-1\}$ and all $i\in\N_0$. The
solution, which depends on the initial value $x_0$, the
perturbation functions $w(\cdot)$ and $v(\cdot)$ and the switching
sequence $\sigma_i$, is denoted by $x_a(n,x_0,v,w)$. 
 
Observe that \eqref{cl2} and \eqref{cla} coincide for 
\[ v(\sigma_i) := \tilde x_{cl}(\sigma_i) - x(\sigma_i), \quad
i=0,1,2,\ldots \] 
for $\tilde x_{cl}$ from Remark \ref{rem:allpredict}.
Similarly, we define
\beq v(n) := \tilde x_{cl}(n) - x_a(n,x_0,v,w), \quad n\ne
\sigma_i . \label{vndef}\eeq 
Note that in contrast to $v(\sigma_i)$,
$i=0,1,2,\ldots$, which are values we are free to choose in
\eqref{cla}, the values in \eqref{vndef} are determined by the
dynamics of the system and the predictor.

The key requirement to establish robustness of the closed loop is now
to assume that the auxiliary system is ISS with respect to $v$ and $w$. 
\begin{assumption}(ISS of the auxiliary system)\\
We assume that the system \eqref{cla} is ISS
with respect to perturbations and 
measurement errors, i.e. there exist $\beta\in\KK\LL$ and
$\gamma_w$, $\gamma_v\in\KK_\infty$ such that 
\[ \|x_a(n,x_0,v,w)\| \le \max\{ \beta(\|x_0\|,n),
\gamma_w(\|w\|_\infty), \gamma_v(\|v\|_\infty)\} \]
holds for each initial value $x_0$, each sequence
$(\sigma_i)_{i\in\N_0}$ satisfying $0\le \sigma_{i+1}-\sigma_i \le m$,
each perturbation function $w\in W$ and each measurement error
 satisfying \eqref{vndef}.
\label{issass}\end{assumption}

Note that for simplicity of exposition we work with a global
 ISS assumption. All subsequent statements can be modified in a
 straightforward way if ISS only holds for sufficiently small
 perturbations and for initial values in a bounded subset of the
 state space. 

 \begin{remark}(ISS of the auxiliary system and MPC) \\
For the MPC setting in Case (iib), Assumption
 \ref{issass} may be optimistic even for $w\equiv 0$ and 
 $v\equiv 0$, since  the controller is computed from an optimization
 over the approximate  prediction \eqref{asys} instead of optimizing
 over the exact  solution of the (nominal) exact system
 \eqref{sys}. Hence,  
 in general, we can only expect stability for the closed loop
 approximate model rather than for the exact one. For
 simplicity of exposition, we work with the simplified
 assumption that the MPC controller stabilizes the
 exact model \eqref{sys}. If desired, this additional error source
 could be rigorously taken into account in the subsequent analysis by  
 formulating Assumption \ref{issass} for an appropriately perturbed
 version of the closed loop approximate system
 \eqref{asys}. Approaches in this direction can be found, for example, in
 \cite{GruN03,Elai07}. These 
 references also show how the needed robustness of the MPC
 controller (and more general optimization based controllers) can be
 obtained from regularity properties of the optimal value function
 which acts here as a Lyapunov function.  
 Alternatively, robustness
 can be ensured by using robust $min-max$ predictive control
 approaches or set based methods, cf. \cite{limon09a}.
 Note that, in general, the gains $\gamma_w$ and
 $\gamma_v$ depend on the prediction horizon and may become larger
 for increasing horizons. 
\label{rem:prederror}
\end{remark}

Now we can establish the  following theorem.
\begin{theorem} (Stability and robustness)\\ Consider a prediction consistent control scheme
 and assume that the controller using this scheme satisfies 
 Assumption \ref{issass}. Then, the trajectories of \eqref{cl2} satisfy\\[-.5cm]
\bea  \|x_{cl}(n,x_0,\sigma_{\cdot},\tau,w)\| \le \max\Big\{ 
\beta(\|x_0\|,n),
    \gamma_w(\|w\|_\infty),\nonumber\\[1ex]
     \gamma_v\Big(\eps(\tau_\infty+\Delta^\sigma_\infty,0) + 
        \eta(\tau_\infty+\Delta^\sigma_\infty,\|w\|_\infty)
        \Big)\Big\} \qquad\eea
 for all $x_0\in X$, $n\in\N_0$ and $w\in W$, 
where $\tau_\infty = \max_{i\in\N_0} \tau(\sigma_i)$ and 
$\Delta^\sigma_\infty = \max_{i\in\N_0}\sigma_{i+1}-\sigma_i$.
\label{allDeltathm}\end{theorem}
{\bf Proof:} for a fixed initial value $x_0$, perturbation $w$ and
switching times $\sigma_i$ we denote the solution of \eqref{cl2} by
$x_{cl}(n)$ and let $\tilde x_{cl}(n)$ be the corresponding predictions  
from Remark \ref{rem:allpredict}.By setting $v(n) = \tilde x_{cl}(n) - x_{cl}(n), \quad n\in\N_0, $
the right hand sides of \eqref{cl2} and \eqref{cla} coincide, and
consequently we achieve $ x_{cl}(n) = x_a(n,x_0,v,w).$
From Assumption \ref{issass}, it follows that in order to prove the
desired inequality we need to show 
\beq \|v\|_\infty \le \eps(\tau_\infty+\Delta^\sigma_\infty,0) + 
        \eta(\tau_\infty+\Delta^\sigma_\infty,\|w\|_\infty). 
\label{desbound}\eeq
For
$n\in\{\sigma_i,\ldots,\sigma_{i+1}-1\}$, it follows from 
\eqref{xpalt} that the predictions satisfy  
$ \tilde x_{cl}(n) = \tilde
x(n,\sigma_i-\tau(\sigma_i),x_{cl}(\sigma_i-\tau(\sigma_i)),\tilde
u).$
Fixing $n$ and $\sigma_i$ and abbreviating
$n_0=\sigma_i-\tau(\sigma_i)$ we obtain that
\beq n-n_0 = n - \sigma_i + \tau(\sigma_i) \le \Delta^\sigma_\infty +
\tau_\infty.\label{nn0bound}\eeq
By using the control input sequence $u$ from \eqref{ucl}, the closed loop
trajectory $x_{cl}$ satisfies $ x_{cl}(n) = x(n,n_0,x_{cl}(n_0),u,w)$.
Since the scheme is prediction consistent, the control sequence $\tilde
u$ used in the prediction coincides with the control sequence $u$ from
\eqref{ucl}. Hence, from Assumption \ref{errass}(i)-(ii), and Equation
\eqref{nn0bound} we have \\[-.7cm] 
\beaa  && \hspace*{-7mm} \|\tilde x_{cl}(n) \!- \!x_{cl}(n)\| \\&=&
\|\tilde
x(n,n_0,x_{cl}(n_0),u)\! -\! x(n,n_0,x_{cl}(n_0),u,w)\| \\
& \le & \| \tilde
x(n,n_0,x_{cl}(n_0),u) - x(n,n_0,x_{cl}(n_0),u,0)\| \\
&& + \es \|x(n,n_0,x_{cl}(n_0),u,0) - x(n,n_0,x_{cl}(n_0),u,w)\|\\
& \le & \eps(n-n_0,\|x_{cl}(n_0) - x_{cl}(n_0)\|) +
\eta(n-n_0,\|w\|_\infty) \\
& \le & \eps(\Delta^\sigma_\infty +
\tau_\infty,0) + \eta(\Delta^\sigma_\infty +
\tau_\infty,\|w\|_\infty). \eeaa
This proves \eqref{desbound} and thus the claim.\endproof

This theorem establishes robustness bounds on the closed loop
system; (16) shows the dependence of the
resulting error with respect to  delays and other factors.

Both $\eps$ and $\eta$ are usually monotonically increasing in their first
argument, cf.\ \eqref{epsex} and \eqref{etaex} and the discussion
after these formulas. Hence, the sensitivity of the
closed loop scheme with respect to the perturbation $w$
crucially depends on the value $\Delta^\sigma_\infty + \tau_\infty$.  

In the scheme described above the $\sigma_i$-sequences are determined
by the network properties: every time the network
is available a new control sequence is sent. Thus, for each $i$ the
difference $(\sigma_{i+1}-\sigma_i)$ is chosen 
as small as possible. Conversely, $(\sigma_{i+1}-\sigma_i)$ and thus
$\Delta^\sigma_\infty$ becomes the larger the 
longer the network is unavailable. The delays
$\tau(\sigma_i)$, on the other hand, are 
determined by the speed of the information transfer: the longer the delay $\tau_s$
from sensor to controller and the longer the anticipated delay
$\tau^{\max}$ from controller to actuator, the larger $\tau_\infty$
becomes. 
\section{Example}\label{sec:example}

We illustrate our result considering the following, simple fourth
order system (two double ``integrators'')
\[\dot x(t) = ( x_2(t),\,  u_1(t),\, x_4(t),\,  u_2(t))^T \]
with state $x=(x_1,\,x_2,\,x_3,\,x_4)^T\in\R^4$ and control
$u=(u_1,\,u_2)^T\in\R^2$. The system is controlled by an MPC
controller with sampling time $T=0.1$ and cost functional  
$ J(x_0, u) = \int\limits_{0}^{NT} l(x(t)) dt +  F(x(NT))$ with
stage cost 
$l(x(t)) = 100 ( x_1^2(t) + x_3^2(t) - 36 )^2 + 0.05 (
x_2(t) + \frac{10 x_3(t)}{\sqrt{x_1^2(t) + x_3^2(t)}} )^2 + 0.05 ( x_4(t) -
\frac{10 x_1(t)}{\sqrt{x_1^2(t) + x_3^2(t)}} )^2$ and terminal cost $F(x)=20
l(x)$, i.e. the closed loop trajectories are supposed
to evolve counterclockwise on the circle $\|(x_1,x_3)^T\|_2=6$ with
constant speed  
$\|(x_2,x_4)^T\|_2=10$. The state and 
control constraints $x_2^2 + x_4^2 \leq 30$ and $u_1^2 + u_2^2 \leq
100$ were imposed, the optimization horizon was $N=10$ and
no stabilizing terminal constraints were used.

The closed loop behavior is simulated with random additive errors in
each state component uniformly distributed in the interval
$[-0.1,0.1]$, where we used the same random sequence for all
simulations. For the simulation 
the prediction consistent scheme from \cite{Grune2009c} is used 
with different delay bounds $\tau^{\max}$. Additionally, communication
failures in the channel from sensor to controller are considered; here
only every third transmission is successful.  This
means that $n_c-n_s = 2$ occurs for every third computation time $n_c$
and consequently $\tau_\infty= \tau^{\max}+2$.  
\begin{figure}[htb]
\centering
\includegraphics[width=0.5\textwidth]{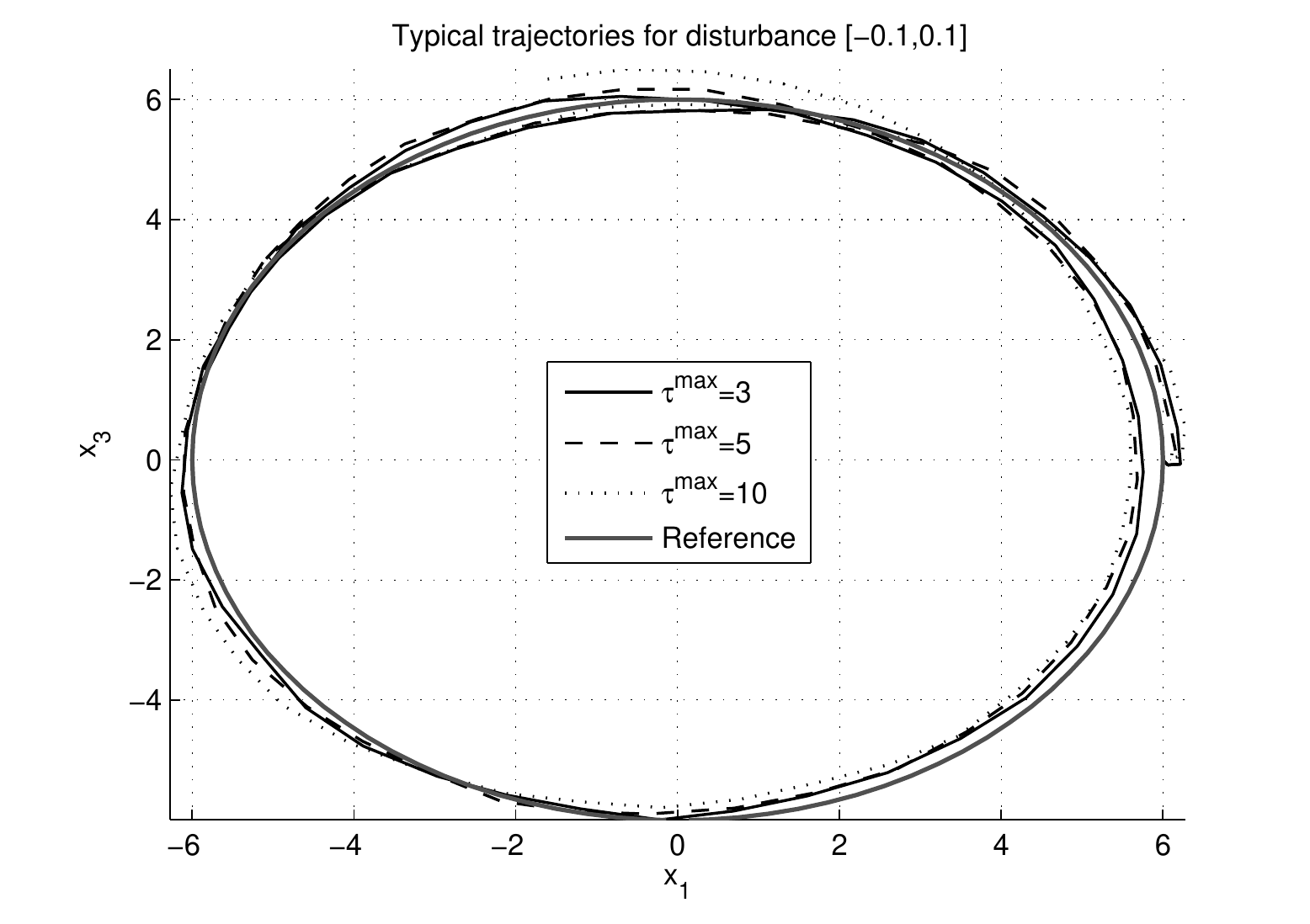}
\caption{Closed loop trajectories in the $(x_1,x_3)$-plane.}\label{fig:x}
\end{figure}

\begin{figure}[htb]
\centering
\includegraphics[width=0.5\textwidth]{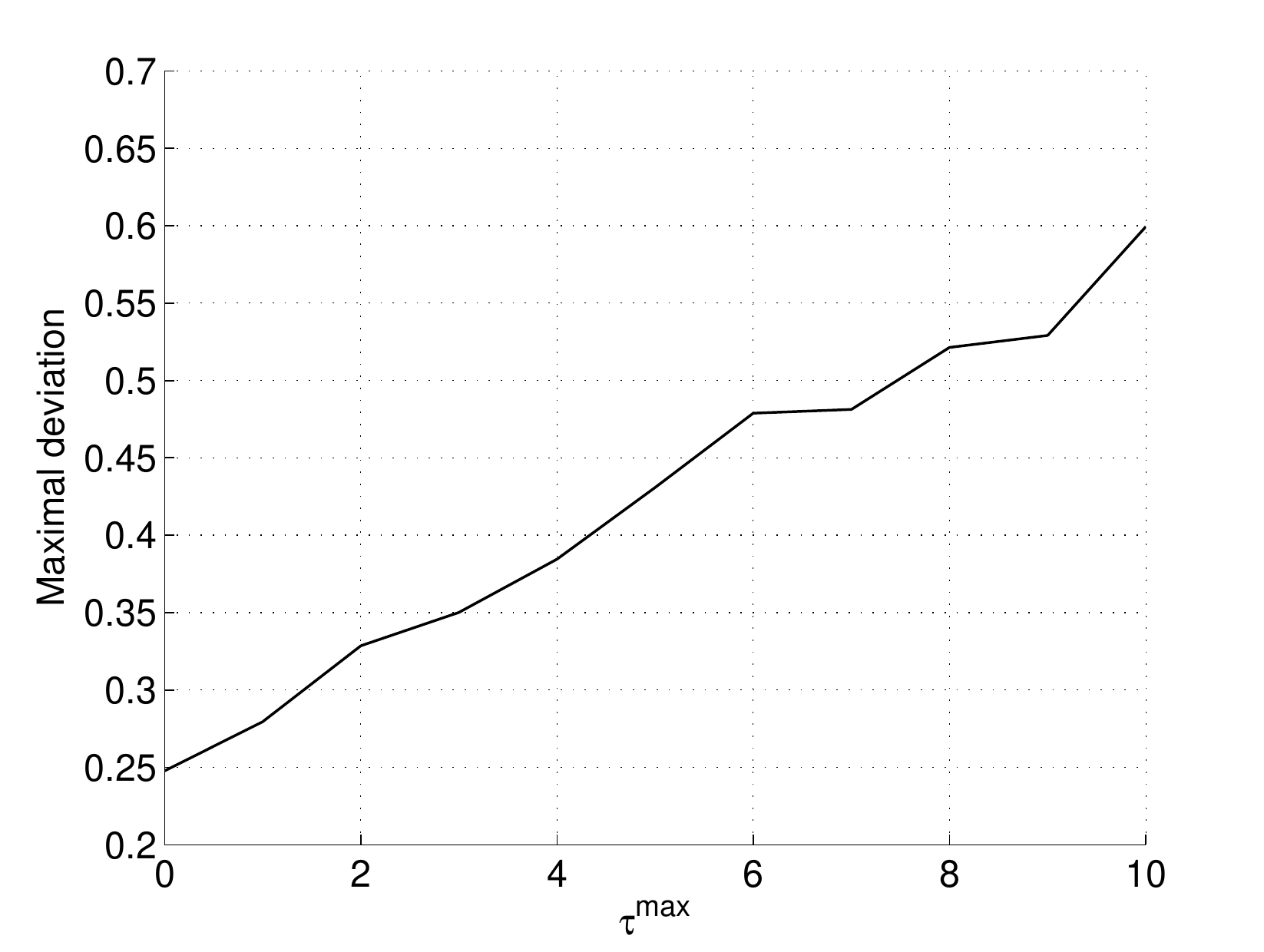}
\caption{Maximal deviation depending on $\tau^{\max}$.}\label{fig:error}
\end{figure}

Fig.\ \ref{fig:x} shows typical trajectories 
and Fig.\ \ref{fig:error} the
maximal deviation 
$\sqrt{(\|(x_1,x_3)^T\|_2-6)^2 + (\|(x_2,x_4)^T\|_2-10)^2}$ depending on
$\tau^{\max}$. 
The results confirm the
robustness of the closed loop as well as the increasing sensitivity against
perturbations for larger $\tau^{\max}$. Note that the deviation 
grows linearly in
$\tau^{\max}$ since the system is open loop stable, cf.\ Remark
\ref{rem:open-loop-stable}.  

\section{Conclusions}\label{sec:conclusions}
Delays and information loss are often unavoidable in todays control
systems. 
They must be accounted for during
the controller design to avoid instability or performance decay. One
way to improve the performance and to guarantee stability is the use
of prediction based compensation schemes. Instead of a single input a
sequence of (predicted) future controls is submitted, buffered and implemented
sequentially at the actuator. The main result of this work is that  
so-called prediction consistent compensation and control schemes, such
as, e.g., the predictive control schemes presented in  
\cite{Findeisen2009} or \cite{Grune2009c}, can posses inherent 
robustness  properties. Precisely, if the nominal closed loop system
without 
delays is input-to-state stable (ISS), then the closed loop system
subject to delays and utilizing a prediction consistent delay
compensation approach is robustly stable. Additionally, we explicitly
analyzed the quantitative influence of the prediction horizon and model
uncertainties on the robustness gains and illustrate the results by a
numerical simulation. The derived results significantly expand the
applicability  of prediction consistent delay compensation approaches,
e.g. based on predictive control solutions, since they establish that
these methods are also applicable to uncertain systems subject to
disturbances.   

\bibliography{references}
\end{document}